\documentclass[11pt]{article}
\PassOptionsToPackage{obeyspaces}{url}

\usepackage{amsfonts,amsmath,amssymb,amsthm}
\usepackage{mathtools}
\usepackage{enumitem}
\usepackage{graphicx}
\usepackage{tikz}
\usetikzlibrary{arrows.meta,calc}
\usepackage{bm}
\usepackage{cancel}
\usepackage{comment}
\usepackage{etoolbox}
\usepackage{xstring}

\usepackage[hidelinks]{hyperref}
\usepackage{bookmark}
\bookmarksetup{open,numbered,depth=5}
\usepackage[nameinlink]{cleveref}

\allowdisplaybreaks[4]

\patchcmd{\thebibliography}
{\leftmargin\labelwidth}
{\leftmargin\labelwidth\addtolength\itemsep{-0.3\baselineskip}}
{}{}

\theoremstyle{plain}
\newtheorem{theorem}{Theorem}[section]

\newtheorem{lemma}[theorem]{Lemma}

\newtheorem{proposition}[theorem]{Proposition}

\theoremstyle{definition}

\newtheorem{problem}[theorem]{Problem}

\crefname{fact}{fact}{facts}
\Crefname{fact}{Fact}{Facts}
\crefname{claim}{claim}{claims}
\Crefname{claim}{Claim}{Claims}
\crefname{observation}{observation}{observations}
\Crefname{observation}{Observation}{Observations}
\crefname{question}{question}{questions}
\Crefname{question}{Question}{Questions}
\crefname{construction}{construction}{constructions}
\Crefname{construction}{Construction}{Constructions}
\crefname{exercise}{exercise}{exercises}
\Crefname{exercise}{Exercise}{Exercises}
\crefname{problem}{problem}{problems}
\Crefname{problem}{Problem}{Problems}

\newcommand*{\arXiv}[1]{\href{https://arxiv.org/abs/#1}{arXiv:#1}}

\definecolor{mred}{RGB}{217,38,38}
\definecolor{mgreen}{RGB}{38,217,38}
\definecolor{mblue}{RGB}{38,38,217}
\definecolor{mcyan}{RGB}{38,128,128}
\definecolor{mmagenta}{RGB}{128,38,128}
\definecolor{myellow}{RGB}{128,128,38}
\definecolor{mgray}{RGB}{98,98,98}

\title{Exact codegree Tur\'an densities of even zycles}
\author{Shuaichao Wang\thanks{\raggedright%
		Faculty of Mathematics and Statistics, Central China Normal University,
		Wuhan 430079, PR China. \texttt{shuaichao@ccnu.edu.cn}.}}
\date{}
\hypersetup{pdftitle={Exact codegree Turan densities of even zycles},pdfauthor={Shuaichao Wang}}

\begin{document}
	\maketitle
	
\begin{abstract}
	We determine the exact codegree Tur\'an density
	of the $3$-uniform zycle $Z_\ell$ for every
	even integer $\ell\ge 8$ with $3\nmid\ell$,
	showing that $\gamma(Z_\ell)=1/3$.
	More generally, we prove that
	$\gamma(Z_\ell)\le 1/3$ for all even $\ell\ge 8$.
	Our results partially answer a question of
	Piga and Sch\"ulke
	[Bull. Lond. Math. Soc. 58 (2026), e70348]
	concerning the sharpness of their upper bound
	of $1/2$.
\end{abstract}

\section{Introduction}\label{sec:introduction}

For a $3$-uniform hypergraph $H$, or $3$-graph for short, the
\emph{codegree} of a pair of vertices is the number of edges
containing that pair. The \emph{minimum codegree} $\delta_2(H)$
is the minimum of these quantities over all pairs of vertices
of $H$. For a fixed $3$-graph $F$, the \emph{codegree Tur\'an
	number} $\operatorname{ex}_2(n,F)$ is the maximum of $\delta_2(H)$
over all $F$-free $3$-graphs $H$ on $n$ vertices. The corresponding
\emph{codegree Tur\'an density} is defined by
\[
\gamma(F)=\lim_{n\to\infty}
\frac{\operatorname{ex}_2(n,F)}{n}.
\]
The existence of this limit was proved by Mubayi and
Zhao~\cite{MZ}. The codegree Tur\'an problem is a local version
of the classical Tur\'an problem, in which the number of edges
is replaced by the minimum codegree.

Determining codegree Tur\'an densities is difficult even for
small $3$-graphs. For the complete $3$-graph $K_4^{(3)}$ on four
vertices, the conjecture of Czygrinow and Nagle~\cite{CN} that
$\gamma(K_4^{(3)})=1/2$ remains open. For $K_4^{(3)-}$, obtained
by deleting one edge from $K_4^{(3)}$, Falgas-Ravry, Pikhurko,
Vaughan and Volec~\cite{FPVV} proved that
$\gamma(K_4^{(3)-})=1/4$ using flag algebras, confirming an
earlier conjecture of Nagle.
Other notable exact results include the codegree
Tur\'an densities of the Fano plane, determined by
Mubayi~\cite{MubayiFano}, and the $3$-graph $F_{3,2}$,
determined by Falgas-Ravry, Marchant, Pikhurko and
Vaughan~\cite{FMPV}.
For an overview of related problems,
see~\cite{BCL2022}.

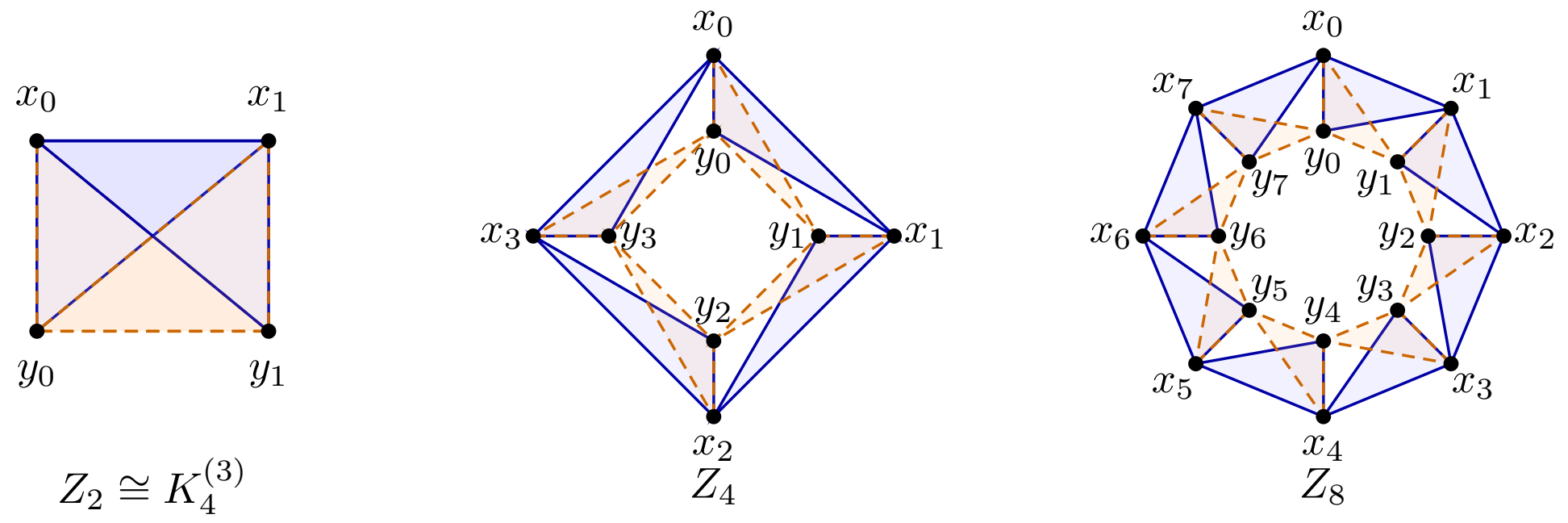
\begin{figure}[htbp]
	\centering
	\begin{tikzpicture}[font=\small,
		vx/.style={circle,fill=black,inner sep=1.25pt},
		ex/.style={draw=blue!65!black,fill=blue!45,
			fill opacity=.12,line width=.65pt},
		ey/.style={draw=orange!80!black,fill=orange!55,
			fill opacity=.12,line width=.65pt,densely dashed}]
		\begin{scope}[xshift=1.7cm]
			\coordinate (x0) at (-.95,.78);
			\coordinate (y0) at (-.95,-.78);
			\coordinate (x1) at (.95,.78);
			\coordinate (y1) at (.95,-.78);
			\foreach \i/\j in {0/1,1/0} {
				\path[ex] (x\i)--(y\i)--(x\j)--cycle;
				\path[ey] (x\i)--(y\i)--(y\j)--cycle;
			}
			\foreach \i in {0,1} {
				\node[vx] at (x\i) {};
				\node[vx] at (y\i) {};
				\node[above=3pt] at (x\i) {$x_\i$};
				\node[below=3pt] at (y\i) {$y_\i$};
			}
			\node at (0,-2.05) {$Z_2\cong K_4^{(3)}$};
		\end{scope}
		\foreach \ell/\shift in {4/6.3,8/11.3} {
			\begin{scope}[xshift=\shift cm]
				\pgfmathtruncatemacro{\last}{\ell-1}
				\foreach \i in {0,...,\last} {
					\coordinate (x\i) at
					({90-360*\i/\ell}:1.48);
					\coordinate (y\i) at
					({90-360*\i/\ell}:.86);
				}
				\foreach \i in {0,...,\last} {
					\pgfmathtruncatemacro{\j}{mod(\i+1,\ell)}
					\path[ex] (x\i)--(y\i)--(x\j)--cycle;
					\path[ey] (x\i)--(y\i)--(y\j)--cycle;
				}
				\foreach \i in {0,...,\last} {
					\node[vx] at (x\i) {};
					\node[vx] at (y\i) {};
					\node at ({90-360*\i/\ell}:1.74)
					{$x_\i$};
					\node at ({90-360*\i/\ell}:.61)
					{$y_\i$};
				}
				\node at (0,-2.05) {$Z_\ell$};
			\end{scope}
		}
	\end{tikzpicture}
	\caption{The zycles $Z_2$, $Z_4$ and $Z_8$.
		Each shaded triangle represents a hyperedge: solid blue
		triangles represent $\{x_i,y_i,x_{i+1}\}$, and dashed
		orange triangles represent $\{x_i,y_i,y_{i+1}\}$,
		with indices taken modulo $\ell$.
		Only the marked points are vertices; crossings do not
		create vertices.}
	\label{fig:zycles}
\end{figure}

In this paper, we study the codegree Tur\'an densities of zycles.
For an integer $\ell\ge2$, the \emph{zycle} $Z_\ell$ is the
$3$-graph with vertex set
$V(Z_\ell)=\{x_i,y_i:i\in\mathbb{Z}_\ell\}$ and edge set
$E(Z_\ell)=\bigl\{\{x_i,y_i,x_{i+1}\},\{x_i,y_i,y_{i+1}\}:
i\in\mathbb{Z}_\ell\bigr\}$,
where the $2\ell$ vertices are distinct and the indices are taken
modulo $\ell$. In particular, $Z_2\cong K_4^{(3)}$, so the
Czygrinow--Nagle conjecture concerns the first member of this
family. See \Cref{fig:zycles} for illustrations of $Z_2$, $Z_4$
and $Z_8$.

Zycles played a central role in the work of Piga and
Sch\"ulke~\cite{PS} on the possible values of codegree Tur\'an
densities. They proved that for every $\varepsilon>0$ there
exists an integer $\ell\ge3$ such that
$0<\gamma(Z_\ell)<\varepsilon$, establishing that zero is an
accumulation point of the set of codegree Tur\'an densities of
$3$-graphs. This contrasts with the classical Tur\'an density
of a $3$-graph, which cannot lie in $(0,2/9)$ by a theorem of
Erd\H{o}s~\cite{Erdos}. Subsequently, Li, Liu, Sch\"ulke and
Sun~\cite{LLSS} used zycles in a recursive construction to show
that $(r-1)/r$ is an accumulation point of the set of codegree
Tur\'an densities of $k$-graphs for all integers $k\ge3$ and
$r\ge1$.

Piga and Sch\"ulke~\cite{PS} also proved that
$\gamma(Z_\ell-e)=0$ for every $\ell\ge3$ and every
$e\in E(Z_\ell)$. The proof of this vanishing-density result
in~\cite{PS} builds on methods developed by
Piga, Sales and Sch\"ulke~\cite{PSS2023}
for tight cycles with one edge removed.
This result was also recovered by Ding,
Lamaison, Liu, Wang and Yang~\cite{DLLWY} through their study
of layered $3$-graphs. Sarkies~\cite{Sarkies} extended it to
$k$-uniform zycles for every $k\ge3$.
Moreover, the conjectured characterization of $3$-graphs with vanishing
codegree Tur\'an density proposed in~\cite{DLLWY} was subsequently
established independently by Yang, Fang and Chen~\cite{YFC}
and King, Piga, Sales and Sch\"ulke~\cite{KPSS}.

For the densities of the zycles themselves, Piga and
Sch\"ulke~\cite[Section~4]{PS} proved that
$\gamma(Z_\ell)\le1/2$ for every $\ell\ge3$ and asked whether
this bound is sharp for some $\ell\ge3$. They also observed
that a cyclic three-part construction gives
$\gamma(Z_\ell)\ge1/3$ whenever $3\nmid\ell$.
Our main result determines the exact density for every even
$\ell\ge8$ not divisible by $3$ and gives an improved upper
bound for all even $\ell\ge8$.
In particular, $\gamma(Z_8)=1/3$. 

\begin{theorem}\label{thm:main}
	For every even integer $\ell\ge8$ with $3\nmid\ell$,
	we have $\gamma(Z_\ell)=1/3$.
	More generally, $\gamma(Z_\ell)\le1/3$ for every
	even integer $\ell\ge8$.
\end{theorem}

\paragraph{Our approach.} 
Our proof is combinatorial. We encode homomorphisms from
zycles as closed walks in an auxiliary digraph on vertex pairs.
Combining averaging arguments with estimates on intersections
of pair neighbourhoods, we find closed walks of every even
length at least eight. Codegree supersaturation then yields
the required copies of zycles.

In \Cref{sec:covers}, we establish two auxiliary lemmas for
digraphs. We then study closed walks in pair digraphs in
\Cref{sec:separation}. The proof of \Cref{thm:main} is given
in \Cref{sec:proof}. We conclude in \Cref{sec:conclusion}
with several open problems.
	
	\section{Auxiliary results on digraphs}\label{sec:covers}
	Walks are directed and may repeat vertices and arcs. For a digraph $D$,
	$N_D^{+j}(v)$ denotes the endpoints of walks of length \emph{exactly}
	$j$ starting at $v$, and $N_D^+(v)=N_D^{+1}(v)$.
	An \emph{exact-time $t$ incoming cover} is a set $S$ such that every
	vertex has a walk of length exactly $t$ ending in $S$.
	We say that a digraph $D$ has the $(m,q)$-common-successor
	property if, for any $m$ vertices (not necessarily distinct),
	some $q$ of them have a common out-neighbour.

	\begin{lemma}\label[lemma]{lem:five}
	Suppose $D$ has $N$ vertices, $\delta^+(D)\ge N/9$,
	and satisfies the \{$(19,7)$, $(12,5)$, $(6,3)$,
	$(3,2)$\}-common-successor properties.
		Then $D$ has an exact-time five incoming cover with at most two vertices.
	\end{lemma}

\begin{proof}
	Choose $v_1,\ldots,v_{19}$ independently and uniformly
	from $V(D)$, and let $S=\{v_1,\ldots,v_{19}\}$.
	Define
$
	U=\{u\in V(D):N_D^+(u)\cap S=\emptyset\}.
$
	Since $\delta^+(D)\ge N/9$, we have
	\[
	\mathbb E|U|
	=\sum_{u\in V(D)}
	\left(1-\frac{|N_D^+(u)|}{N}\right)^{19}
	\le N(8/9)^{19}<N/9.
	\]
	Fix a choice of $S$ for which $|U|<N/9$.
	As every vertex has at least $N/9$ out-neighbours,
	each vertex has an out-neighbour outside $U$,
	which in turn has an out-neighbour in $S$.
	Thus every vertex reaches $S$ in exactly two steps.
	
	We next reduce the number of target vertices in three steps.
	Regard the chosen vertices $v_1,\ldots,v_{19}$ as
	labelled vertices, allowing repetitions.
	By successively applying the
\{$(19,7)$, $(12,5)$, $(6,3)$,
$(3,2)$\}-common-successor properties,
	we can partition these nineteen labels into groups
	of sizes $7,5,3,2,1,1$, such that each group
	of size at least two has a common out-neighbour.
	For each singleton, choose an arbitrary out-neighbour.
	This yields at most six vertices that every vertex
	of $S$ can reach in one step.
	
	Applying the $\{(6,3), (3,2)\}$-common-successor
	properties to these six labels, we similarly obtain
	at most three vertices after one further step.
	Finally, the $(3,2)$-common-successor property
	reduces these three labels to at most two vertices
	in one more step.
	
	Consequently, every vertex of $D$ reaches
	a fixed set of at most two vertices in exactly
	five steps, as required.
\end{proof}

We also need the following weighted covering lemma.
	\begin{lemma}\label[lemma]{lem:weighted}
		Let $J$ be a finite digraph, with loops permitted, and let $\mu$ be a
		probability measure on $V(J)$. Suppose
	$
		\mu(N_J^+(u))>\frac14$ for every $u\in V(J),
		$
		and suppose $J$ satisfies the $\{(2,2), (4,3)\}$-common-successor properties.
		Then $J$ has an exact-time four incoming cover consisting of a single
		vertex. Consequently it contains a closed walk of every integer length
		at least four.
	\end{lemma}

\begin{proof}
	Choose $v_1,\ldots,v_5$ independently according to
	$\mu$, and let $S=\{v_1,\ldots,v_5\}$.
	Define
$
	U=\{u\in V(J):N_J^+(u)\cap S=\emptyset\}.
	$
	Since $\mu(N_J^+(u))>1/4$ for every $u\in V(J)$,
	we have
	\[
	\mathbb E\mu(U)
	=\sum_{u\in V(J)}
	\mu(u)(1-\mu(N_J^+(u)))^5
	\le (3/4)^5<1/4.
	\]
	Fix a choice of $S$ for which $\mu(U)<1/4$.
	For every $u\in V(J)$, there exists
	$w\in N_J^+(u)\setminus U$, since
	$\mu(N_J^+(u))>1/4$.
	As $w\notin U$, it has an out-neighbour in $S$.
	Thus every vertex of $J$ reaches $S$ in
	exactly two steps.
	
	We next reduce the five labels to a single vertex
	in two further steps. Regard $v_1,\ldots,v_5$
	as labelled vertices, allowing repetitions.
	By applying the $(4,3)$-common-successor property
	to any four labels, we find a common out-neighbour
	of three of them. The remaining two labels have
	a common out-neighbour by the
	$(2,2)$-common-successor property.
	This leaves at most two vertices after one step.
	Applying the $(2,2)$-common-successor property
	once more, we obtain a vertex $o$ that every
	vertex of $S$ reaches in exactly two steps.
	Consequently, every vertex of $J$ reaches $o$
	in exactly four steps.
	
	Finally, since every vertex has positive outdegree,
	for any $k\ge 4$, we can take a directed walk
	of length $k-4$ starting at $o$, and then return
	to $o$ in exactly four steps.
	This gives a closed walk of length $k$.
\end{proof}

\section{Closed walks in pair digraphs}\label{sec:separation}

In this section, we study an auxiliary digraph associated
with a 3-graph, whose closed walks encode homomorphisms
from zycles. Our main goal is to establish the existence
of closed walks of all sufficiently large even lengths.

	For a $3$-graph $H$ and a pair $P\in\binom{V(H)}2$, write
$
	N_H(P)=\{v\in V(H)\setminus P:P\cup\{v\}\in E(H)\}.
$
	The \emph{pair digraph} $D(H)$ has vertex set $\binom{V(H)}2$, with
	\[
	P\to Q\quad\Longleftrightarrow\quad Q\subseteq N_H(P).
	\]
	It is loopless, and $N_{D(H)}^+(P)=\binom{N_H(P)}2$.

We first establish a crossing-arc property for
pair neighbourhoods with small intersection.

\begin{lemma}\label[lemma]{lem:cross}
	Let $H$ be an $n$-vertex $3$-graph with
	$\delta_2(H)>n/3+2$. For two pair states
	$p,q$, write $A=N_H(p)$ and $B=N_H(q)$.
	If $|A\cap B|\le 1$, then there exist
	$P\in\binom{A}{2}$ and $Q\in\binom{B}{2}$
	such that $P\to Q$ or $Q\to P$ in $D(H)$.
\end{lemma}

\begin{proof}
	Write $\delta=\delta_2(H)$ and suppose that
	no such arc exists. Let $B'=B\setminus A$,
	$a=|A|$ and $b=|B'|$. Then $A$ and $B'$
	are disjoint, with $a\ge\delta$ and
	$b\ge\delta-1$. By assumption, every pair
	in $A$ has at most one neighbour in $B'$,
	and every pair in $B'$ has at most one
	neighbour in $A$.
	
	We count the codegrees of all pairs with
	one vertex in $A$ and the other in $B'$.
	Edges with two vertices in $A$ contribute
	at most $2\binom{a}{2}$, while those with
	two vertices in $B'$ contribute at most
	$2\binom{b}{2}$. The remaining edges
	contribute at most $ab(n-a-b)$. Hence
	\[
	ab\delta
	\le ab(n-a-b)
	+2\binom{a}{2}+2\binom{b}{2}.
	\]
	Dividing by $ab$, we obtain
	\[
	\begin{aligned}
		\delta
		&\le n-a-b+\frac{a-1}{b}
		+\frac{b-1}{a}\\
		&\le n-2\delta+1
		+\frac{n-2}{\delta-1}\\
		&<n-2\delta+4,
	\end{aligned}
	\]
	where the last inequality follows from
	$\delta>n/3+2$. This implies
	$3\delta<n+4$, a contradiction.
\end{proof}

	We now use the crossing-arc property to find
	even closed walks when two pair states have
	disjoint three-step out-neighbourhoods.

	\begin{lemma}\label[lemma]{lem:separated}
		Let $H$ be an $n$-vertex $3$-graph with
		$\delta_2(H)>n/3+2$, and write $D=D(H)$.
		Suppose that there exist two pair states
		$p,q\in V(D)$ such that
		$
		N_D^{+3}(p)\cap N_D^{+3}(q)=\emptyset.
	$
		Then $D$ contains a closed walk of every
		even length at least eight.
	\end{lemma}

\begin{proof}
	Write $D=D(H)$, $\delta=\delta_2(H)$, and
	\[
	X_i=N_D^{+i}(p),\qquad
	Y_i=N_D^{+i}(q)\quad (i=1,2,3),
	\qquad A=N_H(p),\quad B=N_H(q).
	\]
	Since $D$ has positive minimum outdegree, the assumption
	$X_3\cap Y_3=\emptyset$ implies
	\begin{equation}\label{eq:separated}
		X_i\cap Y_i=\emptyset\qquad (i=1,2,3).
	\end{equation}
	Indeed, any common state in an earlier layer could be
	extended to a common state in the third layer.
	In particular, since $X_1=\binom A2$ and
	$Y_1=\binom B2$, we have $|A\cap B|\le 1$.

We first observe that no three states of $D$
can have pairwise disjoint out-neighbourhoods.
Suppose otherwise that there exist three pair states
$P_1,P_2,P_3$ whose out-neighbourhoods in $D$
are pairwise disjoint. Since
$
N_D^+(P_i)=\binom{N_H(P_i)}{2},
$
we have
$
|N_H(P_i)\cap N_H(P_j)|\le 1$
for $1\le i<j\le 3$.
As $|N_H(P_i)|\ge\delta$ for each $i$,
inclusion--exclusion gives
\[
n\ge\left|\bigcup_{i=1}^3 N_H(P_i)\right|
\ge 3\delta-3,
\]
contrary to $\delta>n/3+2$.

We claim that any two distinct states in $X_1$
have a common out-neighbour.
Indeed, let $P_1,P_2\in X_1$ be distinct and
take any $Q\in Y_1$. Since
\[
N_D^+(P_1),N_D^+(P_2)\subseteq X_2,
\qquad N_D^+(Q)\subseteq Y_2,
\]
and $X_2\cap Y_2=\emptyset$, the
out-neighbourhood of $Q$ is disjoint from
those of $P_1$ and $P_2$.
Thus $P_1$ and $P_2$ must have a common
out-neighbour, as otherwise these three states
would have pairwise disjoint out-neighbourhoods,
a contradiction.
The same holds for any two distinct states in $Y_1$.
	
	Since $|A\cap B|\le 1$, \Cref{lem:cross} guarantees
	an arc between $X_1$ and $Y_1$.
	After interchanging $p$ and $q$ if necessary,
	we may assume that $x\to y$ for some
	$x\in X_1$ and $y\in Y_1$. Thus
	$
	y\in X_2\cap Y_1.
	$

	\medskip
	\noindent\textbf{Claim 1.}
	There exists $z\in X_1\cap Y_2$, and
$
	X_1\cap X_2=Y_1\cap Y_2=\emptyset.
$
	
	\noindent\emph{Proof of Claim 1.}
	We first show that $Y_1\cap Y_2=\emptyset$.
	Otherwise, let $u\in Y_1\cap Y_2$.
	Since $y\in X_2$ and $X_2\cap Y_2=\emptyset$,
	we have $y\ne u$.
	As $y,u\in Y_1$, they have a common
	out-neighbour, which must belong to
	$X_3\cap Y_3$, a contradiction.
	Hence $Y_1\cap Y_2=\emptyset$, or equivalently,
	\begin{equation}\label{eq:B-sparse}
		|N_H(Q)\cap B|\le 1
		\qquad (Q\in Y_1).
	\end{equation}
	
	Now consider the three sets $A$, $B$ and
	$N_H(y)$, each of size at least $\delta$.
	By \eqref{eq:B-sparse},
	$|N_H(y)\cap B|\le 1$, while
	$|A\cap B|\le 1$.
	If $|A\cap N_H(y)|\le 1$, then
	inclusion--exclusion would give
	$n\ge 3\delta-3$, a contradiction.
	Therefore $|A\cap N_H(y)|\ge 2$.
	Choose
	\[
	z\in\binom{A\cap N_H(y)}2.
	\]
	Then $z\in X_1\cap Y_2$.
	
	Finally, suppose that $u\in X_1\cap X_2$.
	Since $z\in Y_2$ and $X_2\cap Y_2=\emptyset$,
	we have $u\ne z$.
	As $u,z\in X_1$, they have a common
	out-neighbour, which would belong to
	$X_3\cap Y_3$, a contradiction.
	Thus $X_1\cap X_2=\emptyset$, or equivalently,
	\begin{equation}\label{eq:A-sparse}
		X_1\cap X_2=\emptyset,\qquad
		|N_H(P)\cap A|\le 1
		\quad (P\in X_1).
	\end{equation}
	\hfill$\diamond$
	
	\medskip
	\noindent\textbf{Claim 2.}
	Every state in $X_2$ has an out-neighbour in $X_1$.
	
	\noindent\emph{Proof of Claim 2.}
	Fix the state $z\in X_1\cap Y_2$ obtained
	in Claim 1, and let $Q\in X_2$.
	Since the out-neighbourhoods of $Q$ and $z$
	lie in the disjoint sets $X_3$ and $Y_3$,
	respectively, we have
	\[
	|N_H(Q)\cap N_H(z)|\le 1.
	\]
	Also, $|N_H(z)\cap A|\le 1$ by
	\eqref{eq:A-sparse}.
	If $|N_H(Q)\cap A|\le 1$, then the three sets
	$N_H(Q)$, $N_H(z)$ and $A$, each of size at
	least $\delta$, would have pairwise intersections
	of size at most one, contradicting
	$n<3\delta-3$.
	Thus $|N_H(Q)\cap A|\ge 2$, which implies
	that $Q$ has an out-neighbour in $X_1$.
	\hfill$\diamond$
	
	\medskip
	By Claim 2, for each $Q\in X_2$, we may choose
	\begin{equation}\label{eq:selection}
		f(Q)\in\binom{N_H(Q)\cap A}{2}
		\subseteq X_1.
	\end{equation}
	In particular, $Q\to f(Q)$.
	
	We now construct an auxiliary digraph $J$ on $X_1$.
	For $P,P'\in X_1$, let $P\to_J P'$ whenever
	there exists $Q\in X_2$ such that
	$P\to_D Q\to_D P'$.
	Note that loops are permitted in $J$.
	Put
	\[
	W=V(H)\setminus A,\qquad w=|W|,
	\qquad S_P=N_H(P)\setminus A
	\quad (P\in X_1).
	\]
	By \eqref{eq:A-sparse} and the minimum
	codegree condition,
	\begin{equation}\label{eq:majority}
		|S_P|\ge\delta-1,\qquad
		w\le n-\delta,\qquad
		\delta-1>\frac{w+1}{2}.
	\end{equation}
	Moreover, $w\ge\delta-1\ge 2$,
	since $S_P\subseteq W$.

	Let $R=X_2\cap\binom W2$, and define a
	probability measure $\mu$ on $X_1$ by
	\[
	\mu(T)=
	\frac{|f^{-1}(T)\cap R|}{|R|}
	\qquad (T\subseteq X_1).
	\]
	Note that $R\ne\emptyset$, since
	$\binom{S_P}{2}\subseteq R$ for every $P\in X_1$.
	
	For every $P\in X_1$ and
	$Q\in\binom{S_P}{2}$, we have
	$P\to_D Q\to_D f(Q)$.
	Thus $f(Q)\in N_J^+(P)$, and hence
	\begin{equation}\label{eq:quarter}
		\mu(N_J^+(P))
		\ge\frac{\binom{|S_P|}{2}}{|R|}
		\ge\frac{\binom{|S_P|}{2}}{\binom w2}
		>\frac14.
	\end{equation}
	The last inequality follows from
	\eqref{eq:majority}.

	It remains to verify the two common-successor properties
	required by \Cref{lem:weighted}.
	For any two states $P_1,P_2\in X_1$, not necessarily
	distinct, we have
	\[
	|S_{P_1}\cap S_{P_2}|
	\ge 2(\delta-1)-w
	\ge 3\delta-n-2>1.
	\]
	Choose a pair
	$Q\in\binom{S_{P_1}\cap S_{P_2}}{2}$.
	Then $f(Q)$ is a common out-neighbour of
	$P_1$ and $P_2$ in $J$.
	
	Next, consider four states $P_1,\ldots,P_4\in X_1$,
	again allowing repetitions.
	For each $v\in W$, let $m_v$ denote the number
	of indices $i\in[4]$ for which $v\in S_{P_i}$.
	Since $\binom{m_v}{3}\ge m_v-2$, we obtain
	\[
	\begin{aligned}
		\sum_{I\in\binom{[4]}{3}}
		\left|\bigcap_{i\in I}S_{P_i}\right|
		&=\sum_{v\in W}\binom{m_v}{3}\\
		&\ge 4(\delta-1)-2w\\
		&\ge 6\delta-2n-4>4.
	\end{aligned}
	\]
	Thus some three of the four sets have at least
	two common vertices. Choosing a pair $Q$
	in their intersection, we obtain a common
	out-neighbour $f(Q)$ of the corresponding
	three states in $J$.
	
	By \Cref{lem:weighted}, $J$ contains a closed
	walk of every integer length $k\ge 4$.
	Replacing each arc of such a walk by its
	corresponding two-step walk in $D$ yields
	a closed walk of length $2k$ in $D$,
	as required.
	
\end{proof}

\section{Proof of the main theorem}\label{sec:proof}

A homomorphism from a 3-graph $F$ to a 3-graph $H$
is a map $\varphi:V(F)\to V(H)$ that sends every
edge of $F$ to an edge of $H$.
It need not be injective on the whole vertex set.
We first establish a homomorphism version of
\Cref{thm:main}.

	\begin{theorem}\label{thm:hom}
		There exists an absolute constant $C_0>0$ such that every
		$n$-vertex $3$-graph $H$ with $\delta_2(H)>n/3+C_0$
		admits a homomorphism $Z_\ell\to H$ for every even integer $\ell\ge8$.
	\end{theorem}
	
	Fix $C_0>2$ sufficiently large that
	\[
	mC_0\ge\binom{m}{q}
	\qquad
	\text{for }(m,q)\in
	\{(3,2),(6,3),(12,5),(19,7)\}.
	\]
	We first verify that the pair digraph $D(H)$
	satisfies the assumptions of \Cref{lem:five}.

	\begin{lemma}\label[lemma]{lem:common}
		If $\delta_2(H)>n/3+C_0$, then $D(H)$ satisfies
		the $\{(19,7),(12,5),(6,3),(3,2)\}$-common-successor properties, and has
		minimum outdegree greater than $N/9$,
		where $N=\binom{n}{2}$.
	\end{lemma}
	\begin{proof}
		Write $\delta=\delta_2(H)$ and $N=\binom n2$. Then
		\[
		\delta^+(D(H))\ge\binom\delta2
		>\frac{(n/3+1)(n/3)}2>\frac N9.
		\]
		
		For $m$ pair states $P_1,\ldots,P_m$, not necessarily
		distinct, let
		\[
		m_v=|\{i\in[m]:v\in N_H(P_i)\}|
		\qquad (v\in V(H)).
		\]
		Since $\binom{m_v}{q}\ge m_v-(q-1)$, double counting gives
		\[
		\begin{aligned}
			\sum_{I\in\binom{[m]}q}
			\left|\bigcap_{i\in I}N_H(P_i)\right|
			&=\sum_{v\in V(H)}\binom{m_v}{q}\\
			&\ge\sum_{v\in V(H)}m_v-(q-1)n\\
			&=\sum_{i=1}^m|N_H(P_i)|-(q-1)n\\
			&\ge m\delta-(q-1)n.
		\end{aligned}
		\]
		
		For each $(m,q)\in\{(3,2),(6,3),(12,5),(19,7)\}$, we have
		$m/3\ge q-1$. By the choice of $C_0$,
		\[
		m\delta-(q-1)n
		>\left(\frac m3-q+1\right)n+mC_0
		\ge mC_0\ge\binom mq.
		\]
		Thus some intersection contains at least two vertices, whose pair
		is a common successor of the corresponding $q$ labels.
	\end{proof}
	
	\begin{proof}[Proof of \Cref{thm:hom}]
		By \Cref{lem:common,lem:five}, the pair digraph $D=D(H)$ has an
		exact-time five incoming cover $S$ with $|S|\le2$.
		If $S=\{p\}$, then $p$ has a closed walk of every integer length at
		least five: first walk the required excess length, then use the cover
		to return in five steps.
		
		Otherwise write $S=\{p,q\}$. If
		$
		N_D^{+3}(p)\cap N_D^{+3}(q)\ne\varnothing,
	$
		choose a common endpoint $z$. Every vertex reaches $p$ or $q$ in
		exactly five steps, and then $z$ in exactly three further steps.
		Thus $\{z\}$ is an exact-time eight incoming cover and $D$ has
		closed walks of every integer length at least eight.
		If the two three-step sets are disjoint, \Cref{lem:separated} gives
		closed walks of all even lengths at least eight.
		
		In every case, for even $\ell\ge8$, take such a closed walk and map
		the two vertices of the $i$th distinguished pair of $Z_\ell$ to the
		two vertices of its $i$th pair state. Each arc guarantees the two
		required hyperedges. Adjacent states are disjoint, and nonconsecutive
		states may overlap, as permitted in a homomorphism. Hence $Z_\ell\to H$.
	\end{proof}

	To complete the proof of \Cref{thm:main},
	we pass from homomorphisms to copies of zycles
	using codegree supersaturation.
	For a nonempty finite family $\mathcal F$ of $3$-graphs, let
	$\operatorname{ex}_2(n,\mathcal F)$ be the largest minimum codegree
	of an $n$-vertex $3$-graph containing no member of $\mathcal F$,
	and put
	\[
	\gamma(\mathcal F)
	=\lim_{n\to\infty}\frac{\operatorname{ex}_2(n,\mathcal F)}{n}.
	\]
	This limit exists by~\cite{MZ}.
	For a $3$-graph $F$, let $F(t)$ denote its uniform $t$-blow-up.
	We shall use the following codegree supersaturation result
	of Piga and Sch\"ulke~\cite[Proposition~2.4]{PS}.
	
	\begin{proposition}\label[proposition]{prop:supersaturation}
		Let $\mathcal F$ be a nonempty finite family of
		$3$-graphs. For every $\varepsilon>0$ and
		every positive integer $t$, there exists $n_0$
		such that every $n$-vertex $3$-graph $H$ with
		$n\ge n_0$ and
		\[
		\delta_2(H)\ge(\gamma(\mathcal F)+\varepsilon)n
		\]
		contains a copy of $F(t)$ for some $F\in\mathcal F$.
	\end{proposition}
	
	\begin{proof}[Proof of \Cref{thm:main}]
		Fix an even integer $\ell\ge 8$, and let
		$\mathcal F_\ell$ be the family of all
		homomorphic images of $Z_\ell$, up to isomorphism.
		Note that $\mathcal F_\ell$ is finite.
		By \Cref{thm:hom}, every $\mathcal F_\ell$-free
		$n$-vertex $3$-graph has minimum codegree
		at most $n/3+C_0$. Consequently,
		\[
		\gamma(\mathcal F_\ell)\le\frac13.
		\]
		Moreover, every homomorphism from $Z_\ell$
		to some $F\in\mathcal F_\ell$ lifts to an
		injective embedding of $Z_\ell$ into
		the $(2\ell)$-blow-up $F(2\ell)$, by assigning
		distinct clones to vertices with the same image.
		
		Now fix $\varepsilon>0$. By
		\Cref{prop:supersaturation}, every sufficiently
		large $n$-vertex $3$-graph $H$ with
		\[
		\delta_2(H)\ge(1/3+\varepsilon)n
		\]
		contains a copy of $F(2\ell)$ for some
		$F\in\mathcal F_\ell$, and hence a copy
		of $Z_\ell$. Therefore,
		$\gamma(Z_\ell)\le 1/3$.
		
		For the matching lower bound when $3\nmid\ell$,
		consider the cyclic three-part construction
		of Piga and Sch\"ulke~\cite[Section~4]{PS}.
		Partition the vertex set into three parts
		$A_0,A_1,A_2$ of sizes differing by at most one,
		and take all triples with two vertices in
		$A_i$ and one in $A_{i+1}$, where indices
		are taken modulo three.
		The resulting $3$-graph $H$ satisfies
		$
		\delta_2(H)\ge\lfloor n/3\rfloor-1.
		$
		
		Observe that the neighbourhood of every pair
		in $H$ is contained in a single part.
		Consequently, every state on a closed walk
		in $D(H)$ is a pair contained in some $A_i$.
		Moreover, every arc from a pair in
		$\binom{A_i}{2}$ leads to a pair in
		$\binom{A_{i+1}}{2}$.
		Hence every closed walk in $D(H)$ has
		length divisible by three.
		
		Since every homomorphism from $Z_\ell$ to
		$H$ induces a closed walk of length $\ell$
		in $D(H)$, the construction is $Z_\ell$-free
		whenever $3\nmid\ell$.
		It follows that
		$\gamma(Z_\ell)\ge 1/3$, completing the proof.
	
	\end{proof}

	\section{Concluding remarks}\label{sec:conclusion}
	We have shown that $\gamma(Z_\ell)\le 1/3$ for every
	even integer $\ell\ge 8$, with equality whenever
	$3\nmid\ell$. Several questions remain open.
	
	Since $Z_2\cong K_4^{(3)}$, determining $\gamma(Z_2)$
	is the classical problem of Czygrinow and Nagle~\cite{CN},
	who conjectured that $\gamma(Z_2)=1/2$.
	
For the remaining small even zycles $Z_4$ and $Z_6$, the results
	of Piga and Sch\"ulke~\cite[Section~4]{PS} give
	\[
	\frac13\le\gamma(Z_4)\le\frac12,
	\qquad
	\frac14\le\gamma(Z_6)\le\frac12.
	\]
	It is natural to ask whether our upper bound extends
	to these two cases.
	
	\begin{problem}\label{prob:short-even}
		Determine $\gamma(Z_4)$ and $\gamma(Z_6)$.
		In particular, is $\gamma(Z_4)=1/3$,
		and is $\gamma(Z_6)\le1/3$?
	\end{problem}
	Our approach does not directly extend to odd zycles,
	since \Cref{lem:separated} guarantees only even closed
	walks. In particular, determining $\gamma(Z_3)$
	was already posed by Piga and
	Sch\"ulke~\cite[Problem~4.2]{PS}.
	Moreover, when $3\mid\ell$, the cyclic three-part
	construction does not provide a matching lower bound.
	This leads to the following problem.
	
	\begin{problem}
		Determine $\gamma(Z_\ell)$ for odd $\ell\ge 3$
		and for even $\ell\ge 12$ divisible by three.
	\end{problem}

\end{document}